\documentclass[12pt]{article}

\usepackage{amsmath,amssymb,url}
\usepackage{tikz,color,pgf}
\usepackage{longtable} 

\newtheorem{theorem}{Theorem}[section]
\newtheorem{corollary}[theorem]{Corollary}
\newtheorem{lemma}[theorem]{Lemma}

\newtheorem{conjecture}[theorem]{Conjecture}
\newtheorem{claim}[theorem]{Claim}
\newtheorem{obs}[theorem]{Observation}

\newcommand{\proof}{\noindent{\bf Proof.\ }}
\newcommand{\qed}{\hfill $\square$ \bigskip}
\newcommand{\w}{{\rm s}}

\newcommand{\smallqed}{{\tiny ($\Box$)}}
\newcommand{\sa}{{\rm s}}
\newcommand{\ggd}{\gamma_g}

\begin{document}

\title{$1/2$-conjectures on the domination game and claw-free graphs
}

\author{Csilla Bujt\'as$^{a}$\thanks{Email: \texttt{csilla.bujtas@fmf.uni-lj.si}} 
\and Vesna Ir\v si\v c$^{a,b}$\thanks{Email: \texttt{vesna.irsic@fmf.uni-lj.si}}
\and Sandi Klav\v zar $^{a,b,c}$\thanks{Email: \texttt{sandi.klavzar@fmf.uni-lj.si}}
}
\maketitle

\begin{center}
$^a$ Faculty of Mathematics and Physics, University of Ljubljana, Slovenia\\
\medskip

$^b$ Institute of Mathematics, Physics and Mechanics, Ljubljana, Slovenia\\
\medskip

$^c$ Faculty of Natural Sciences and Mathematics, University of Maribor, Slovenia\\
\medskip
\end{center}

\begin{abstract}
	Let $\gamma_g(G)$ be the game domination number of a graph $G$. Rall conjectured that if $G$ is a traceable graph, then $\gamma_g(G) \le \left\lceil \frac{1}{2}n(G)\right\rceil$. Our main result verifies the conjecture over the class of line graphs. Moreover, in this paper we put forward the conjecture that if $\delta(G) \geq 2$, then $\gamma_g(G) \leq \left\lceil \frac{1}{2}n(G) \right\rceil$. We show that both conjectures hold true for claw-free cubic graphs. We further prove the upper bound  $\ggd(G) \le \left\lceil \frac{11}{20} \, n(G) \right\rceil$ over the class of claw-free graphs of minimum degree at least $2$. Computer experiments supporting the new conjecture and sharpness examples are also presented. 
\end{abstract}


\noindent
{\bf AMS Subj.\ Class.\ (2020)}: 05C57, 05C69, 05C76

\section{Introduction}
\label{sec:intro}

Papers~\cite{bresar-2010, kinnersley-2013} are the milestones for the success of the domination game. The first of them introduced the game, proved some fundamental results, and described the imagination strategy which turned out to be one of the main proof techniques for the game. The second paper delivered the Continuation Principle which is a tool used to derive many important results, and the $3/5$-conjecture which asserts that if $G$ is an isolate-free graph, then $\gamma_g(G)\le \frac{3}{5}n(G)$, where $n(G)$ denotes the order of $G$.  These two papers caused the area to flourish---well above fifty papers have already been written. We emphasize papers~\cite{dorbec-2015, klavzar-2019, kosmrlj-2017, nadjafi-2016, ruk-2019, xu-2018} and investigations of different variants of the game~\cite{borowiecki-2019, bresar-2019, bujtas-2016, duchene-2020, henning-2015, henning-2017, jiang-2019}. 

A lot of research has been done on the $3/5$-conjecture, the first results being reported in~\cite{kinnersley-2013}. In~\cite{bresar-2013} the conjecture was verified for all trees on at most $20$ vertices and all trees which meet the conjectured bound with equality were listed. A large family of extremal trees was subsequently presented in~\cite{henning-2017b} and it was conjectured that there are no other such trees. Proofs of many results on the $3/5$-conjecture as a key tool use the discharging method, the method being initiated in~\cite{bujtas-2015a}. In 2016, Henning and Kinnersley~\cite{henning-2016} proved a fundamental result that the conjecture holds for all graphs with minimum degree at least $2$. For bigger minimum degrees, better upper bounds were earlier proved in~\cite{bujtas-2015b}. In particular, we know that $\gamma_g(G)< 0.4803\, n(G)$ holds whenever $\delta(G) \ge 5$. The best general upper bound was recently published in~\cite{bujtas-2020a}: if $G$ is an isolate-free graph, then $\gamma_g(G)\le \frac{5}{8}n(G)$. Further, Rall conjectured that the $3/5$-conjecture can be improved for traceable graphs as follows. 

\begin{conjecture}
\label{conj:Rall}
If $G$ is a traceable graph, then $\gamma_g(G) \le \left\lceil \frac{1}{2}n(G)\right\rceil$.
\end{conjecture}

Rall's conjecture was published for the first time in~\cite{james-2019}, see also~\cite{bujtas-2020b}. Now, knowing that if $\delta(G)\ge 2$, then $\gamma_g(G) \leq \frac{3}{5}n(G)$ holds, and on the other hand not being aware of minimum degree $2$ graphs $G$ for which $\gamma_g(G) > \left\lceil \frac{1}{2}n(G) \right\rceil$ would hold, we pose in this paper the following: 

\begin{conjecture}
	\label{conj:1/2-for-all}
	If $\delta(G) \geq 2$, then $\gamma_g(G) \leq \left\lceil \frac{1}{2}n(G) \right\rceil$.
\end{conjecture}

In this paper we provide different supports for both conjectures and proceed as follows. In the next section we first give definitions and concepts needed, and then present a set up for our proofs, as well as the intuition behind it. In Section~\ref{sec:claw-free} we prove that both conjectures hold true for claw-free cubic graphs. In Section~\ref{sec:Ham-line-graphs} we turn our attention to Hamiltonian line graphs, another rich subclass of claw-free graphs. We first consider the edge domination game which was, to the best of our knowledge, considered by now only in the unpublished paper~\cite{Ta-2014}. We prove upper bounds on the game edge domination number for graphs with edge dominating circuits or with edge dominating trails. As a consequence we deduce that Conjecture~\ref{conj:Rall} is valid on the class of line graphs. In Section~\ref{sec:claw-free-mindeg-2} we consider claw-free graphs of minimum degree $2$ and prove for them a weakened version of Conjecture~\ref{conj:1/2-for-all}, that is, $\ggd(G) \le \left\lceil \frac{11}{20} \, n(G) \right\rceil$ holds for claw-free graphs $G$ with $\delta(G) \ge 2$. In the concluding section we present our computer experiments and a variety of examples that demonstrate that Conjecture~\ref{conj:1/2-for-all} is sharp, of course provided it holds true.  

\section{Preliminaries}
\label{sec:preliminaries}

The domination game on a graph $G$ is played by Dominator and Staller. If Dominator starts the game, we speak of a D-game, otherwise it is an S-game. During the game, the players alternately select vertices that are not dominated by the set of previously selected vertices. When no such vertex is available the game is over. Dominator's goal is to end the game as soon as possible, Staller's goal is the opposite. The unique number of moves played in the D-game when both players play optimally is the {\em game domination number} $\gamma_g(G)$. In the S-game the corresponding graph invariant is the {\em Staller-start game domination number} $\gamma_g'(G)$ of $G$. 

If $v$ is a vertex of a graph $G$, then the \emph{open neighborhood} $N_G(v)$ is the set of neighbors of $v$, while the \emph{closed neighborhood} is $N_G[v] = N_G(v)\cup \{v\}$. Vertices $u$ and $v$ are \emph{true twins} in $G$, if $N_G[u]=N_G[v]$. The closed neighborhood of a set $S\subseteq V(G)$ of vertices is $N_G[S]=\bigcup_{v \in S}N_G[v]$. If $G$ will be clear from the context, we may omit the subscript $G$ in this notation. 

A \emph{claw} is a star $K_{1,3}$. A graph is called \emph{claw-free} if it does not contain an induced subgraph isomorphic to a claw. A  \emph{diamond} is a $K_4-e$. The two vertices of a diamond that are of degree $3$ (and incident with two triangles) are called the \emph{central vertices of the diamond}. Note that the two central vertices of a diamond are true twins in every cubic graph. 

Denoting by $e_1, \dots, e_m$ the edges of a graph $G$, the \emph{line graph} $L(G)$ of $G$ contains one vertex $v_i$ for every  $e_i$, $i\in [m]$.  Two vertices $v_i$ and $v_j$ are adjacent in $L(G)$ if $e_i$ and $e_j$ share a vertex in $G$. 
 By definition, the vertices in $L(G)$ represent the edges of $G$ so that every edge dominating set in $G$ corresponds to a (vertex) dominating set in $L(G)$. 

\subsection{On the proof methods} \label{sec:proof-meth}

We first introduce the terminology that will be used throughout the proofs of our main theorems. At the end of this subsection, we provide an informal hint about the ideas and working methods that may result a potential function which can be successfully applied for proving upper bounds on graph invariants.
\medskip

Given a graph $G$ and a set $D \subseteq V(G)$, the \emph{(vertex) residual graph} $G^D$ is the graph $G$ together with the following color assigment to its vertices:

\begin{itemize}
	\item A vertex $v$ is \emph{white} if $v \notin N[D]$.
	\item  A vertex $v$ is \emph{blue} if $v \in N[D]$ and $N[v] \not\subseteq N[D]$.
	\item A vertex $v$ is \emph{red} if $N[v] \subseteq N[D]$.
\end{itemize}
Moreover, in $G^D$, the sets of white, blue and red vertices, respectively, are denoted by $W$, $B$, and $R$. In the domination game, $D$ is interpreted as the set of vertices which have been played until a point in the game. Then, $v$ is undominated if and only if $v \in W$, while $u$ is a legal move in $G^D$ if and only if $v \in W \cup B$. The game starts with $W=V(G)$ in $G^\emptyset$, and ends when $R=V(G)$. The latter equivalently means $W=\emptyset$ and $N[D]=V(G)$.

For a vertex $v$ in the residual graph $G^D$, the \emph{white-degree} (or shortly, \emph{W-degree}) of $v$ is the number of its white neighbors, that is, $\deg_W(v)= |N[v] \cap W|$. By definition, $\deg_W(v)=0$ if $v$ is red, and $\deg_W(v)\ge 1$ if $v$ is blue. The maximum W-degree over the white vertices of $G^D$ is denoted by $\Delta_W(W)$, while $G^D[W]$ stands for the subgraph induced by the white vertices in $G^D$. 

In the proofs we introduce a \emph{weight function} $f^D: V(G) \rightarrow \mathbb{R}$ for each $G^D$, where the weight of a vertex $v$ depends on its color and, maybe, on some further properties. The \emph{weight of the residual graph}  is defined as $f(G^D)= \sum_{v\in V(G)}f^D(v)$.  It is also called the  \emph{potential function} over the residual graphs.
If $G^D$ is fixed and $v \in V(G)$, then $\sa(v)$ denotes the decrease in the potential function when $v$ is played in $G^D$ (for unplayable vertices $\sa(v)=0$). Formally, $\sa(v)=f(G^D)-f(G^{D\cup\{v\}})$. We specify $f(G^D)$ in the proofs such that $\sa(v)\ge 0$ holds for every $v \in V(G)$.

 In Section~\ref{sec:Ham-line-graphs}, we use an analogous terminology for edges instead of vertices. Two (different) edges $e$ and $f$ are \emph{adjacent} in a simple graph $G$, if they share a vertex that is, if $e\cap f \neq \emptyset$. We say that an edge $e\in E(G)$ \emph{dominates} itself and the adjacent edges. This set of edges is the \emph{closed neighborhood of $e$} defined formally as $N_G[e]=\{f : f \in E(G) \mbox{ and } e\cap f \neq \emptyset\}$. The \emph{open neighborhood} of $e$ is $N_G(e)=N_G[e]\setminus \{e\}$.
 A set $D$ of edges is an \emph{edge dominating set} in $G$ if  
 $$N_G[D]=\bigcup_{e\in D}N_G[e]=E(G).$$ 
 
 For a graph $G$ and a set $D \subseteq E(G)$, the \emph{(edge) residual graph} $G^D$ is the graph $G$ together with the following color assignment. 
 An \emph{edge} $e$ is \emph{white} if $e \notin N[D]$; $e$ is \emph{blue} if $e \in N[D]$ and $N[e] \not\subseteq N[D]$; and $e$ is \emph{red} if $N[e] \subseteq N[D]$ holds.
  For a fixed (edge) residual graph $G^D$, the sets of white, blue and red edges, respectively, are denoted by $W$, $B$, and $R$. 
  An edge $e$ is undominated if and only if $e \in W$ and it can be played in the next turn if and only if $e \in W \cup B$. 
  
  We also introduce an \emph{(edge) weight function} $f^D: E(G) \rightarrow \mathbb{R}$ for each $D \subseteq E(G)$ and the \emph{weight of the residual graph}  is defined as $f(G^D)= \sum_{e\in E(G)}f^D(e)$.  If $G^D$ is fixed and $e \in E(G)$, then let $\sa(e)=f(G^D)-f(G^{D\cup\{e\}})$.

  In the proofs, we consider the vertex or edge domination games with specified  potential functions $f$ such that $f(G^{D})=0$ holds when the game is over. Our aim is to prove that Dominator has a strategy which ensures that the average decrease of $f(G^D)$ in a move of the game exceeds a given constant $c$. This implies that, assuming that the players follow optimal strategies, the number of moves is at most $\lfloor f(G^\emptyset)/c \rfloor$.
  \medskip
  
The formal description of the proofs starts with the definition of a potential function which needs the preceding specification of the weights. This preliminary work on the problem does not appear in the manuscript, and therefore, it might seem a bit mysterious. In fact, the first step is designing an algorithm (or a strategy for Dominator), finding the key points in it, and, by representing the weights of the vertices with variables, constructing a set of linear inequalities. These constraints typically express the aim that Dominator's move, together with Staller's next move, decreases the sum of the weights by at least $1$ on average. The next step is routine. We want to minimize the potential function's value at the beginning of the game while all the mentioned inequalities are satisfied. Usually, it suffices to solve an LP-problem with a linear objective function. 
 	
 	However, the first (and the second, the third) attempt does not typically result in the desired minimum. In this case, we may modify the algorithm and replace some sharp inequalities with weaker constraints. We may also introduce a more subtle weight assignment by making differences between vertices of the same color. Another option is to discharge the weights at a point in the game, or to supplement the potential function by terms depending on some general properties or subgraphs of the residual graph. After including several ideas, we hopefully obtain a potential function that is appropriate for the proof.

\section{Claw-free cubic graphs}
\label{sec:claw-free}

Our main aim in this section is proving Conjecture~\ref{conj:Rall} and  Conjecture~\ref{conj:1/2-for-all} for the class of claw-free cubic graphs.
 
 An edge $e$ is called \emph{triangle edge} if it is incident with a triangle; otherwise, $e$ is a \emph{flat edge}. It is easy to see that in a claw-free cubic graph a vertex $v$ is either contained in a $K_4$-component or $v$ is a central vertex of a diamond or, in the third case, $v$ is incident to exactly one flat edge and one triangle.

\begin{theorem} \label{thm:claw-free-cubic}
	If $G$ is a claw-free cubic graph, then $\ggd(G) \le \frac{n(G)}{2}$.
\end{theorem}

\proof Throughout the proof, we assume that $G$ is a claw-free cubic graph on $n=n(G)$ vertices and $D \subseteq V(G)$ does not dominate $G$. We define a potential function $f$ on the possible residual graphs such that $f$ assigns $3n$ and $0$, respectively, to the residual graph at the beginning and at the end of the game. Then, we describe a strategy for Dominator which ensures that the average decrease of $f$ in a move of the game will be at least $6$ when the game ends. This will directly imply that the number of moves in the game is at most $\frac{3n-0}{6}=\frac{n}{2}$.
\medskip

In a residual graph $G^D$, a diamond subgraph is called K-subgraph if its  central vertices are white and the remaining two vertices are blue. Observe that, in a K-subgraph $S$, only the blue vertices are adjacent to vertices outside $S$. By the definition of blue vertices, all these neighbors outside $S$ must be red as they belong to $D$. This implies that K-subgraphs are pairwise vertex disjoint and there is no edge between them.  

For a residual graph $G^D$, we define the potential function
$$f(G^D)= 3|W|+|B|$$
which can also be interpreted as a weight assignment $f^D: V(G) \rightarrow \{0,1,3\}$. In this assignment, every white (that is, undominated) vertex gets a weight of $3$, while blue and red vertices are assigned to $1$ and $0$, respectively. We will also use the notation $\sa(v)$ for the nonnegative difference $f(G^D)-f(G^{D\cup\{v\}})$.

\medskip
The following claim relates to the case when $v$ is a legal move in the domination game that is, if $v \in W\cup B$. We will use this result when estimate the decrease $\sa(v)$ for a move $v$ of Staller.
\begin{claim} \label{cl:St}
	If $v$ is a legal move on $G^D$ in the domination game, then $\sa(v) \ge 3$.
\end{claim}
\proof If $v \in W$, then $v$ is recolored red and already this change ensures $\sa(v) \ge 3$. If $v \in B$, then $v$ becomes red and since, by definition of a legal move, $v$ has at least one white neighbor in $G^D$ that becomes blue or red in $G^{D\cup \{v\}}$, the decrease $\sa(w)$ is at least $(1-0)+(3-1)=3$. This establishes the claim.
\smallqed
\medskip

In some cases we will need a stricter result on the move of Staller. In the continuation, $G^D[W]$ denotes the subgraph induced by the white vertices of $G^D$.

\begin{claim} \label{cl:St2}
	If $v$ is a move which dominates a vertex $u$ from a $P_1$ component of $G^D[W]$, or from a cycle component $C_i$ with $i \ge 4$ of $G^D[W]$, or from a K-subgraph of $G^D$, then $\sa(v) \ge 5$.
\end{claim}
\proof First suppose that $u$ forms a $P_1$ component in $G^D[W]$. Then, the neighbors $u_1$, $u_2$, and $u_3$ of $u$ are blue. Note that $u$ cannot be a central vertex of a diamond or a vertex from a $K_4$ component as in these cases the true twin neighbor(s) of $u$ would be white as well. Therefore, $u$ is incident with a triangle, say it is $uu_1u_2$, and with a flat edge $uu_3$ in $G$. By definition, each blue vertex has a red neighbor from $D$. Then, the triangle incident with $u_3$ has a red vertex and the only white neighbor of $u_3$ is $u$. The same is true for $u_1$ and $u_2$ since their neighbors outside $uu_1u_2$ must be red. Then, an arbitrary move which dominates $u$ turns $u$, $u_1$, $u_2$, and $u_3$ to be red. As follows, we have $\sa(v) \ge 3+3\cdot 1=6$ in the first case. 

Now, suppose that $u$ belongs to a cycle component $C_i$ of $G^D[W]$ such that $i \ge 4$. This white cycle cannot contain central vertices from a diamond subgraph and therefore, every vertex $x$ from $C_i$ is incident with one flat edge and with one triangle in $G$. Since $C_i$ is not a triangle and it is an induced subgraph in $G^D[W]$, flat and triangle edges alternate along the white cycle $C_i$. Thus, if the played vertex $v$ dominates $u$ and $v \in B$, then there is a  $vuu'$ triangle with  $u,u' \in V(C_i) \subseteq W$. As $u$ and $u'$ are blue in $G^{D\cup \{v\}}$,   $\sa(v) \ge 2\cdot (3-1)+1=5$ holds. In the other case, the played vertex $v$ is white. To dominate $u$, vertex $v$ must belong to  $V(C_i)$. With this move, $v$ becomes red and its neighbors from $V(C_i)$ become blue. We conclude $\sa(v) \ge 3 +2\cdot (3-1) =7$.

Finally, we suppose that the white vertex $u$ is a central vertex of a K-subgraph. When $u$ is dominated, its true twin $u'$ also becomes dominated. Since $N[u]\setminus\{u,u'\} \subseteq B$ and $N[u']\setminus\{u,u'\} \subseteq B$ hold in $G^D$, both $u$ and $u'$ turn into red after the move $v$. This shows that $\sa(v) \ge 2\cdot 3 = 6$.
\smallqed
\medskip

In the main part of the proof we argue that the following property is true for any residual graph $G^D$:
\begin{itemize}
	\item[$(Q)$] If Dominator plays first in $G^D$, he can ensure that at least one of the following statements will be satisfied:
	\begin{itemize}
		\item[$(Q1)$] Dominator finishes the game with the first move and decreases $f(G^D)$ by at least $6$;
		\item[$(Q2)$] The first two moves on $G^D$ decrease $f(G^D)$ by at least $12$;
		\item[$(Q3)$] Dominator finishes the game with the third move and these three moves decrease $f(G^D)$ by at least $18$;
		\item[$(Q4)$] The first four moves on $G^D$ decrease $f(G^D)$ by at least $24$.
	\end{itemize}
\end{itemize}
\begin{claim} \label{cl:2}
	If $D \subset V(G)$ is not a dominating set in $G$, then $G^D$ satisfies property $(Q)$.
\end{claim}
\proof We consider all possibilities concerning the subgraphs of $G^D$ and prove that $(Q)$ is true in every case. By our assumption, $D$ is not a dominating set and hence, there is an undominated vertex; i.e.\ $W\neq \emptyset$. 
\paragraph{(C1)} There exists a white vertex $v$ with $3$ white neighbors.\\
In this case, Dominator plays $v$ that results in the following changes: $v$ becomes red and its neighbors become blue or red. We obtain $\sa(v) \ge 3+ 3(3-1)=9$. If the game finishes with this move, then $(Q1)$ is satisfied. Otherwise, by Claim~\ref{cl:St}, together with the next move of Staller, the potential function is decreased by at least $12$ in two moves and $(Q2)$ holds. 
\medskip

From now on we assume that (C1) is not true for $G^D$. In other words,  $G^D[W]$ has maximum degree at most $2$ and, consequently, it consists of path and cycle components. In particular, each $K_4$ component of $G$ (if exists) is completely red in $G^D$.

\paragraph{(C2)} (C1) is not true and there exists a triangle induced by three white vertices.\\ 
Let $v_1v_2v_3$ be such a triangle in $G^D$. Since (C1) is not true, the neighbor of $v_i$ outside the triangle must be blue, for $i \in [3]$. Thus, if Dominator plays $v_1$, all the three vertices of the triangle are recolored red and we have $\sa(v_1) \ge 3 \cdot 3 =9$. Together with Claim~\ref{cl:St} this prove that $G^D$ satisfies $(Q1)$ or $(Q2)$. 

\paragraph{(C3)} Neither (C1) nor (C2) is true and there is a path component $P_j$ with $j \ge 3$  in  $G^D[W]$.\\ 
Let $P_j=v_1v_2\dots v_j$ be a path component in $G^D[W]$. Since (C1) and (C2) do not hold, $G^D[W]$ is triangle-free and $P_j$ does not contain two edges from the same triangle. Since the two central vertices of a diamond subgraph are true twins, they are of the same color and hence, as $G^D[W]$ is triangle-free, at most two of the four vertices of a diamond  may be white. This proves that $P_j$ does not contain two consecutive triangle edges. As every vertex of $G$ is incident with at most one flat edge, $P_j$ cannot contain two consecutive flat edges either. Consequently, flat and triangle edges alternate along this path. Now, suppose that Dominator plays $v_2$. With this move both $v_1$ and $v_2$ become red since before this move $v_2$ was the only undominated neighbor of $v_1$. Vertex $v_3$ turns to be blue or red. Moreover, one of $v_1v_2$ and $v_2v_3$ is a triangle edge. The third vertex $u$ of this triangle was a blue vertex and had two white neighbors in $G^D$. Hence, $u$ becomes red with the move $v_2$. This yields $\sa(v_2) \ge 2 \cdot 3 + 2 +1 =9$ and we get, by Claim~\ref{cl:St}, that $G^D$ satisfies $(Q1)$ or $(Q2)$.  

\paragraph{(C4)} None of (C1)--(C3) is true and there is a $P_2$ component in  $G^D[W]$ which does not belong to a K-subgraph.\\
First suppose that the $P_2$ component $uv$ is a flat edge in $G$. Then, $u$ is incident with a triangle $uu_1u_2$ where $u_i$ is blue and its third neighbor (outside the triangle $uu_1u_2$) must be red, for $i\in [2]$. The situation is similar for $v$ and for the triangle $vv_1v_2$. Since the graph is cubic and $u_i$ has a red neighbor, $u_i$ is different from $v_j$, and they cannot be adjacent vertices. This is true for every pair $(i,j) \in [2]\times [2]$. Then, the move $u$ of Dominator makes all the six vertices red and $\sa(u) \ge 2\cdot 3 + 4\cdot 1=10$. 

Now, suppose that the $P_2$ component $uv$ is an edge of the triangle $uvw$ where $w$ must be blue. Let $u'$ and $v'$ denote  the other blue neighbors of $u$ and $v$, respectively. Since it is not a K-subgraph, $u'\neq v'$. Then, $uu'$ is a flat edge and $u'$ is dominated via a triangle edge. We may conclude that $u$ is the only white neighbor of $u'$ and that the same is true for $v$ and $v'$.  Remark that $N[w] \cap W=\{u,v\}$ also holds as the third neighbor of $w$ (that is the vertex in $N(w) \setminus \{u,v\}$) has to be red. If Dominator plays $u$, all the five vertices $u$, $v$, $w$, $u'$, $v'$ become red and we have $\sa(u) \ge 2 \cdot 3 +3 \cdot 1 =9$. This ensures that $G^D$ satisfies $(Q)$.

\paragraph{(C5)} None of (C1)--(C4) is true and there exists a K-subgraph in $G^D$.\\
In this case, Dominator may play any vertex $v$ from the K-subgraph, all the four vertices are recolored red and $\sa(v) \ge 2\cdot 3 + 2\cdot 1 = 8$. If the game is over with this move, then $(Q1)$ is satisfied. Otherwise, observe that the vertices of this K-subgraph $S$ were adjacent only to two red vertices from $G^D - S$. Then, it remains true in $G^{D \cup \{v\}}$ that none of (C1)-(C4) is satisfied. In particular, each legal move of Staller in $G^{D \cup \{v\}}$ satisfies the condition in Claim~\ref{cl:St2} and hence, it decreases the potential function by at least $5$. That is, independently of Staller's choice, the two moves together decreases $f(G^D)$ by at least $13$.

\paragraph{(C6)} None of (C1)--(C5) is true and there exists a cycle component in  $G^D[W]$.\\
Since (C2) is not satisfied, this cycle component $C_j$ is of length $j \ge 4$. As it was shown in the proof of Claim~\ref{cl:St2}, flat and triangle edges alternate along this white cycle. It also follows that $j$ is even. Let $C_j=v_1v_2\dots v_j$ such that $v_1v_2$, $v_3v_4$, \dots, $v_{j-1}v_j$ are the triangle edges. $N[V(C_j)]$ contains $j/2$ blue vertices which we denote by $u_2, u_4,\dots,u_j$ such that the triangles incident to the white cycle are $v_1v_2u_2$, $v_3v_4u_4$,\dots,$v_{j-1}v_ju_j$. Each $u_i \in \{u_2, u_4, \dots, u_j\}$ is blue and, by definition, it has a red neighbor from $D$. Consequently, we have  $N[u_i] \cap W=\{v_{i-1}, v_i\}$ for each even $i$, $2 \le i \le j$. Suppose that Dominator plays $v_2$ from $G^D$. Then, $v_2$ and $u_2$ become red, while $v_1$ and $v_3$ become blue. We have $\sa(v_2)= 3 +1 + 2\cdot 2=8$. With this move, the $C_j$ component is replaced with a $P_{j-3}$ component in $G[W]$, but the other components remain untouched. If Staller replies by dominating a vertex from a white component different from $P_{j-3}$, Claim~\ref{cl:St2} implies that the decrease in the potential function is at least $5$. Together with the move of Dominator, $f(G^D)$ is decreased by at least $13$ and $(Q2)$ is satisfied. The same is true if $j=4$ and Staller dominates a vertex from the remaining part of the white cycle $C_4$. Indeed, in this case $P_{j-3}$ is a $P_1$ and we can refer to Claim~\ref{cl:St2} again.

Suppose that $j \ge 6$ and Staller dominates a white vertex from $P'=v_4, \dots, v_j$.  If she plays a white vertex $v_s$ from $P'$,  it becomes red and its white neighbor (which exists as $j \ge 6$) becomes blue or red. Then, $\sa(v_s)\ge 3+2=5$. If Staller plays $u_i$ with $i \ge 6$, as this blue vertex still has two white neighbors, namely $v_{i-1}$ and $v_i$, we have $\sa(u_i)\ge 1+ 2 \cdot 2=5$.  If she plays $v_3$ or $u_4$, the result  is the same: $v_3$ and $u_4$ become red, $v_4$ becomes blue. This gives $\sa(v_3)=\sa(u_4)= 2\cdot 1 + 2=4$. In all these cases, $(Q2)$ is satisfied.

The last possible case is when Staller plays the blue vertex $v_1$ and dominates $v_j$. Then, $v_1$ turns into red and $v_j$ turns into blue but there are no further changes in colors. So, we have only $\sa(v_1)=1+2=3$. But in this situation Dominator can reply by playing $v_{j-2}$. If $j \ge 8$, the white vertices $v_{j-2}$, $v_{j-1}$ and the blue vertices $u_j$, $u_{j-2}$ become red, while  $v_{j-3}$ becomes blue. Then, we have $\sa(v_{j-2})= 2\cdot 3+2\cdot 1 + 2=10$. Together with the next move of Staller which decreases the potential function by at least $3$ according to Claim~\ref{cl:St}, $f(G^D)$ drops by at least $8+3+10+3=24$ and $(Q4)$ is fulfilled. If $j=6$, the move $v_{j-2}=v_4$ of Dominator turns the white vertices $v_{5}$ and  $v_{4}$ into red and the blue vertices $u_6$, $u_{4}$ and $v_6$, $v_3$ also become red. We conclude $\sa(v_4)= 2 \cdot 3+4\cdot 1=10$ again. If this move finishes the game, $(Q3)$ is satisfied, otherwise together with Staller's next move $(Q4)$ holds.

\paragraph{(C7)} None of (C1)--(C6) is true and there exists a $P_1$ component in  $G^D[W]$.\\
If $v$ forms a $P_1$ component in  $G^D[W]$, all its neighbors are blue. Note that $v$ cannot be a central vertex of a diamond subgraph of $G$. If $uvw$ is the triangle containing $v$, then $u$ is dominated from a vertex which is outside of this triangle. Therefore, $N[u]\cap W=\{v\}$ and similarly, $N[w]\cap W=\{v\}$ holds. If $N(v)\setminus \{u,w\}=\{x\}$, then   $vx$ is a flat edge and $x$ is dominated by a vertex of the incident triangle $xx'x''$. Consequently, $x$ cannot have a white neighbor different from $v$. This means that each move which dominates $v$ turns $v$, $u$, $w$, and $x$ into red. For Dominator's move, this gives $\sa(v) = 3+3\cdot 1=6$. If the game finishes with this move, $G^D$ satisfies $(Q1)$. Otherwise, Staller replies in a residual graph where (C7) still holds and her move $y$ gives $\sa(y) =6$ again. We conclude that $(Q)$ is satisfied.
\medskip

If none of (C1)--(C7) holds for the residual graph $G^D$, then $W=\emptyset$ and $D$ is a dominating set in $G$. This finishes the proof of Claim~\ref{cl:2}. 
\smallqed  
\medskip

Starting with $G^\emptyset$ and repeatedly applying Claim~\ref{cl:2} for the residual graphs until the game ends, we get that Dominator may ensure that the average decrease in the potential function in a move is at least $6$. Since $f(G^\emptyset)=3n$ and we have $f(G^D)=0$ when the game finishes, there are at most $3n/6$ moves. We conclude $\gamma_g(G) \le n/2$.
\qed

If Staller starts the domination game on a claw-free cubic graph $G$, her first move $v$ turns the white vertex $v$ into red and turns the three white neighbors into blue or red. Hence, the first move gives $\sa(v) \ge 3+ 3\cdot 2= 9$. Starting with $G^{\{v\}}$ and with the next move of Dominator, we may apply Claim~\ref{cl:2} repeatedly again. The conclusion is that 
$$\gamma_g'(G) \le \frac{3n-9}{6}+1= \frac{n-1}{2}\,.$$

\begin{corollary}
	If $G$ is a claw-free cubic graph, then $\gamma_g'(G) \le \frac{n(G)-1}{2}$.
\end{corollary}

We have thus proved that Conjectures~\ref{conj:Rall} and~\ref{conj:1/2-for-all} hold true for claw-free cubic graphs. With respect to the first  conjecture we argue in the rest of this section that the class of Hamiltonian (and consequently traceable) claw-free graphs is quite rich, so that Theorem~\ref{thm:claw-free-cubic} gives a large, non-trivial class of graphs for which Conjecture~\ref{conj:Rall} holds true. 

Recall that a diamond is defined as a $K_4-e$ that is equivalently a $C_4+e$. Let us now introduce the notation $C_6^+$ for the simple graph obtained from two vertex disjoint triangles by adding two vertex disjoint edges to it.

Suppose that $G$ is a claw-free cubic graph and it is also $(K_4,C_4+e, C_6^+)$-free. Then, every vertex is incident with exactly one triangle and one flat edge. Moreover, any two triangles are vertex disjoint and there is at most one edge between them. These properties ensure that every triangle is adjacent to exactly three other triangles in $G$. Now, define $G^*$ such that each vertex from $V(G^*)$ represents a triangle from $G$ and two vertices are adjacent in $G^*$ if and only if there is an edge between the corresponding triangles in $G$. Observe that $G^*$ is a cubic graph and that every cubic graph $F$ can be obtained this way from a cubic graph $G$ which is $(K_{1,3},K_4,C_4+e, C_6^+)$-free. The latter claim is true even if $F$ contains $K_{1,3},K_4,C_4+e, C_6^+$ subgraphs. To see this, just replace every vertex of $F$ by a triangle and, for every edge $uv \in E(F)$, make adjacent two vertices from the corresponding triangles $t(u)$ and $t(v)$ such that the maximum degree does not  exceed $3$ in the constructed graph.
	
Now consider a claw-free cubic graph $G$ which is also $(K_4,C_4+e, C_6^+)$-free. If $G$ is Hamiltonian, then every triangle must be induced by three consecutive vertices of the Hamiltonian cycle. This directly yields a Hamiltonian cycle for $G^*$ and vice versa, if we are given a Hamiltonian cycle in the cubic graph $F=G^*$, this directly defines a Hamiltonian cycle in $G$. It is proved in~\cite{RW} that, among  $n$-vertex cubic graphs, the proportion of the  Hamiltonian ones tends to $1$ as $n \rightarrow \infty$. Hence the class of the Hamiltonian cubic graphs is quite rich, the same is then true for the Hamiltonian claw-free cubic graphs. We note that the class of the non-Hamiltonian cubic graphs and hence, that of the non-Hamiltonian claw-free cubic graphs are quite rich as well. This fact is demonstrated by NP-completeness of the problem of deciding whether a 2-connected cubic bipartite planar graph admits a Hamiltonian cycle~\cite{AT}.

 \section{Hamiltonian line graphs} \label{sec:Ham-line-graphs}
 
In this section we first revitalize the edge domination game which was inroduced in 2014 in the manuscript~\cite{Ta-2014}, but received no attention afterwards. Then we prove upper bounds on the game edge domination number for graphs with edge dominating circuits or trails. As a consequence we then obtain the main result of this section which asserts that Conjecture~\ref{conj:Rall} is valid on the class of line graphs.

\subsection{Edge domination game}
  
 An open (resp.\ closed) \emph{trail} is an open (resp.\ closed) walk in a graph with no repeated edge. As we consider only simple graphs, a trail will be represented by the corresponding sequence of (not necessarily different) vertices.  A closed trail is also called \emph{circuit}.  	
 An \emph{edge dominating circuit} is a circuit  $v_0, \dots, v_\ell,v_0$  such that every $e \in E(G)$ is dominated by an edge $v_iv_{i+1}$, that is  $e\cap \{v_0, \dots ,v_\ell\} \neq \emptyset$ holds for every $e\in E(G)$. Similarly, an \emph{edge dominating trail} is a trail $v_0, \dots, v_\ell$, the vertex set of which forms a vertex cover in $G$.
 \medskip
 
 The \emph{edge domination game}, first discussed in~\cite{Ta-2014}, is defined analogously to the domination game. Two players, Dominator and Staller, take turns choosing an edge of $G$. Playing an edge $e$ is legal if $N_G[e]$ contains at least one edge not dominated by the edges previously chosen. The game finishes when the set of the played edges becomes an edge dominating set in $G$. Dominator wants to minimize the length of the game (i.e., the number of played edges) while Staller wants to maximize it. The \emph{game edge domination number $\gamma_{e,g}(G)$} (\emph{S-game edge domination number $\gamma'_{e,g}(G)$}) is the length of the game when Dominator (Staller) starts the game and both players follow an optimal strategy according to their goals.
 
By the definition of the line graph, every edge domination game played on $G$ naturally corresponds to a domination game on $L(G)$, and it also holds in the other way around, therefore the following fact holds. 

 \begin{obs}[\cite{Ta-2014}] 
 	\label{obs:1}
 	For every graph $G$ and its line graph $F=L(G)$ we have $\gamma_g(F)= \gamma_{e,g}(G)$ and $\gamma_g'(F)=\gamma_{e,g}'(G)$. 	
\end{obs}

 \subsection{Graphs with edge dominating circuits}
 
We introduce some new concepts before proving our main result concerning the edge domination game on graphs which admit edge dominating circuits. 
 
 For a graph $G$ and a set $D \subseteq E(G)$, we will refer to the (edge) residual graph $G^D$ and related notations as introduced in Section~\ref{sec:proof-meth}. Furthermore, we also assign colors to the vertices based on the colors of the incident edges.  A  vertex $v$ is \emph{white} if it is incident to at least one white edge, otherwise $v$ is \emph{red}. 
 
  Consider an edge $e=uv$ in $G^D$. By definition, if $e$ is white then $u$ and $v$ are white vertices; if $e$ is red then so are $u$ and $v$. If $e$ is blue, then it is adjacent to a white edge and therefore, at least one of $u$ and $v$ is white. Moreover, as $e \in N[D]$, there is an adjacent edge $f \in D$. As follows, the common end of $e$ and $f$ is a red vertex. This gives the following property.
 \begin{obs} \label{obs:2}
 For every $D \subseteq E(G)$, each blue edge of $G^D$ is incident to one white and one red vertex.
 \end{obs}
 
 If we fix an edge dominating circuit $C=v_0\dots v_\ell v_0$ in a graph $G$, the edges of the circuit are called \emph{$C$-edges} while the remaining ones are the \emph{outer edges} in $G$. 

\begin{lemma} \label{lem:2+}
 	If $G$ admits an edge dominating circuit $C$, then every edge $e \in E(G)$ has at least two adjacent $C$-edges. In particular, if $e$ is white in a residual graph $G^D$, then it is adjacent to at least two $C$-edges from $W \cup B$.
\end{lemma}

\proof If $e=v_iv_{i+1}$ is an edge from the circuit $C=v_0\dots v_\ell v_0$, then $v_{i-1}v_i$ and $v_{i+1}v_{i+2}$ are $C$-edges adjacent to $e$.  Since $C$ is an edge dominating circuit, every outer edge $e$ is incident to at least one vertex from $C$. If $v_i$ is such a vertex for $e$, then $e$ is adjacent with the $C$-edges $v_{i-1}v_i$ and $v_{i}v_{i+1}$.
For the second part of the statement, recall that both ends of a white edge are white vertices and hence, there cannot be adjacent red edges. 
\qed

 \begin{theorem} \label{thm:edge-dom-game}
 If $G$ is a graph with an edge dominating circuit, then $$\gamma_{e,g}(G)\le \left\lceil \frac{m(G)}{2}\right\rceil$$
and, if $G$ is not a cycle of length $n \equiv 1 \mod 4$, then 
 $$\gamma_{e,g}(G)\le \left\lfloor \frac{m(G)}{2}\right\rfloor$$
 also holds under the given conditions.
 \end{theorem}
 \proof First suppose that $G$ is a cycle $C_n$. Then its line graph is also isomorphic to $C_n$ and $m(G)=n(L(G))=n$ holds. By Observation~\ref{obs:1}, $\gamma_{e,g}(G)=\gamma_g(L(G))$. Thus,  $\gamma_{e,g}(G)=\gamma_g(C_n)\le \lceil n/2 \rceil$ holds and moreover, $\gamma_{e,g}(C_n)\le n/2$ is also true whenever $n \not\equiv 1 \mod 4$. (We have used here the known formula for $\gamma_g(C_n)$, see~\cite{kosmrlj-2017}.) This verifies the statement for cycles.
 
  From now on, we assume that $G$ is not a cycle but admits an edge dominating circuit $C=v_0\dots v_\ell v_0$ and, therefore, $\Delta(G) \ge 3$. For every $D\subseteq E(G)$, we define a weight function $f^D$ on the edges of $G^D$. 

\begin{center}
	\begin{tabular}{l|c} 
		Type of  $e$	&   $f^D(e)$  \\
		\hline 	
		White edge &  $2$ \\
	Blue $C$-edge & $1$\\
	Blue outer edge & $0$\\
		Red edge &  $0$
	\end{tabular}
\end{center}
The weight of $G^D$ is defined as the sum of the weights assigned to its edges. Formally,
$f(G^D)= \sum_{e\in E(G)}f^D(e)$ is the potential funcion that we consider in the proof.
Further, as introduced earlier, $\w(e)$ stands for the difference $f(G^D)-f(G^{D\cup \{e\}})$. Remark that $w(e)$ is always positive if $e$ is a playable edge in $G^D$.
As for Dominator's moves, we will assume that he always chooses an edge to play  which maximizes the decrease $\w(e)$ in the current residual graph. Our aim is proving that this greedy strategy of Dominator ensures that in every two consecutive moves $d_i$ and $s_i$, 
$f(G^D)$ decreases by at least $8$. We realize this plan by proving a series of claims.

\begin{claim} \label{claim:1}
If a white edge $e$ is played in $G^D$, then every adjacent blue edge becomes red and every adjacent white edge becomes blue or red. If a blue edge $e=uv$ is played and $u$ is a white vertex, then every white edge incident with $u$ becomes blue or red and every blue edge incident with $u$ becomes red. 
\end{claim}
\proof If an edge $e=uv$ is played, then all adjacent edges become dominated (i.e., blue or red) and therefore, $u$ and $v$ will be red in $G^{D\cup \{e\}}$. Suppose first that  $e$ is white in $G^D$ and that it is adjacent to $f \in B$. We may assume, without loss of generality, that $f=uu'$. Then, $u$ is the white end of $f$ and, by Observation~\ref{obs:2}, the other end $u'$ is red. As both $u$ and $u'$ are red in $G^{D\cup \{e\}}$, the edge $f$ becomes red too. The situation is similar if $e$ is blue in $G^D$ and $f$ is adjacent to $e$ in its white end vertex $u$. \smallqed

 \begin{claim} \label{claim:2} 
 Each move of Staller decreases $f(G^D)$ by at least $2$.
 \end{claim}
 \proof Suppose that Staller plays $e$ in $G^D$. Since $e$ is playable, $e \in W\cup B$ in $G^D$ and, with this move, $e$ becomes red. Hence, if $e\in W$ in $G^D$, its weight reduces from $2$ to $0$ that proves $\w(e)\ge 2$. If $e$ is blue and adjacent to an outer white edge $e'$, then $e'$ becomes dominated in $G^{D\cup \{e\}}$. No matter whether $e'$ is blue or red in $G^{D\cup \{e\}}$, its weight will be $0$ there. This fact itself ensures $\w(e)\ge 2$. In the remaining case, Staller plays a blue edge $e$ that is adjacent to a white edge $e_1$ from the circuit. This implies that the white end $u$ of $e$ is incident to the circuit. Let a $C$-edge which is incident to $u$ and different from $e_1$ be denoted by $e_2$. (Note that $e_2= e$ is possible.) When Staller plays $e$, the white  $e_1$ becomes blue or red and this change decreases the weight by at least $1$. If $e_2$ is white then it also turns blue or red and $f^D(e_2)$ drops by at least 1. If $e_2$ is a blue edge in $G^D$ then its white end is $u$ and, by Claim~\ref{claim:1},  $e_2$ will be red in  $G^{D\cup \{e\}}$. In this case, again, $f^D(e_2)$ decreases by $1$ as $e_2$ is a $C$-edge. Therefore, the inequality $\w(e) \ge 1+1=2$ also holds if  $e$ is blue. \smallqed

\begin{claim} \label{claim:3}
  If $G^D$ contains two adjacent white edges, then Dominator can play an edge $e$ such that  $f(G^D)$ decreases by at least $6$.
\end{claim}
\proof Consider the first move $d_1$ of Dominator in $G^\emptyset$. By our assumption, $G$ is not a cycle. If there is an outer edge $e$, let $e_1$   be a $C$-edge  adjacent to $e$ and set $d_1=e_1$. In $G^{D\cup \{e_1\}}$, the edge $e_1$ is red that decreses the weight by $2$. In addition, by Lemma~\ref{lem:2+}, there are at least two $C$-edges adjacent to $e_1$. These edges become blue (or red) that results a further decrease of at least $1+1=2$ in $f(G^\emptyset)$. As $e$ is an outer edge and becomes blue or red, its weight decreases by $2$. This proves $\w(e) \ge 6$. 

In the other case, there are no outer edges and $C$ is not a cycle. Hence, there is a vertex $v_j$ that is incident to at least four $C$-edges, denote them by $e_1,e_2,e_3,e_4$. Each of these four edges is adjacent with an edge wich is not incident to $v_j$. Let $f_1$ be such a neighbor of $e_1$. If Dominator plays $e_1$, it will be red and each of the edges $f_1,e_2,e_3, e_4$ becomes blue or red. This yields $\w(e_1)\ge 2+ 4 \cdot 1=6$ as required.
\medskip

In the later turns of Dominator, $D$ is not empty. Since every edge of $G$ is incident to at least one vertex from $C$, we have some red vertices on the circuit. While $G^D$ contains white edges, there must be white vertices on $C$. Therefore, between consecutive white and red vertices of $C$, we have some blue edges on the circuit. 

First we suppose that the circuit $C$ contains a white edge $v_jv_{j+1}$ such that $v_j$ is incident to more than two $C$-edges. Then, there are at least four $C$-edges being incident to $v_j$ so that each of them belongs to $W\cup B$. Choose four such $C$-edges and denote them by $e=v_jv_{j+1}$, $e_1$, $e_2$ and $e_3$. Observe that $e$ is adjacent with the $C$-edge $f=v_{j+1}v_{j+2}$ and $f\neq e_i$ for every $i \in [3]$. Suppose that Dominator plays $e$. The edge $e$ becomes red that decreases the weight by $2$. For every $i$, if $e_i$ is white in $G^D$, it is recolored blue or red that reduces the weight by at least $1$; if $e_i$ is blue than, by Claim~\ref{cl:St2}, it becomes red and the weight drops by $1$. The situation is the same for the $C$-edge $f$. We may infer that $\w(e) \ge 2+ 4\cdot 1=6$.   From now on we may assume that each white vertex $v_j$  ($j \in [\ell]\cup \{0\}$) is incident to exactly two $C$-edges. Consider the following three cases on the types of adjacent white edges in $G$.
\begin{itemize}
\item Suppose that there exist two adjacent white outer edges $e_1$ and $e_2$. When $e_1$ is played, it is recolored red and $e_2$ becomes blue or red. These changes decrease the weight by $2\cdot 2$. By Lemma~\ref{lem:2+}, there exist two different $C$-edges $f_1$ and $f_2$ that are adjacent to $e_1$. Then, $f_i \in B\cup W$ in $G^D$ for $i \in [2]$. If $f_i\in W$, it becomes blue or red in $G^{D\cup \{e_1\}}$ and therefore, its weight decreases by at least $1$. If $f_i \in B$ then, by Claim~\ref{claim:1}, it becomes red and the weight decreases by $1$. We may infer that $\w(e_1)\ge 2\cdot 2 + 2\cdot 1=6$.
\item Suppose that there exist two adjacent white edges $e_1$ and $e_2$ such that $e_1$ is a $C$-edge and $e_2$ is an outer edge. When $e_1$ is played, $e_1$ becomes red and the outer edge $e_2$ becomes blue or red. This  decreses  the weight by $2\cdot 2$. Further, by Lemma~\ref{lem:2+}, $e_1$ is adjacent to two $C$-edges $f_1$ and $f_2$ both are from $W\cup B$. It does not matter whether $f_i$ is white or blue, its weight decreases by at least $1$ when $e_1$ is played. This proves $\w(e_1) \ge  2\cdot 2 + 2\cdot 1=6$ for the second case.
\item Otherwise, every two adjacent white edges are both $C$-edges.  Recall that, since $D\neq \emptyset$, there are both white and red vertices among $v_0, \dots, v_\ell$. Then, as there exist two adjacent white edges on $C$, we can find a vertex $v_j$ such that the edge $v_{j-1}v_j$ is blue, and both $v_jv_{j+1}$ and $v_{j+1}v_{j+2}$ are white. Therefore, the vertex  $v_{j-1}$ is red while $v_{j}$,  $v_{j+1}$ and  $v_{j+2}$ are white.  By our present condition, $v_j$ is incident to exactly two $C$-edges and it is not  incident to any outer white edges. Thus, $v_jv_{j+1}$ is the only white edge being incident to $v_j$. Consequently, if the edge $v_jv_{j+1}$ becomes dominated, then the vertex $v_j$ turns red and, additionally, the edge $v_{j-1}v_j$ is recolored red as both of its vertices become red.
 We conclude that when $e=v_{j+1}v_{j+2}$ is played,  all the three vertices $v_j,v_{j+1}, v_{j+2}$ and the three edges $v_{j-1}v_{j}$, $v_{j}v_{j+1}$, $e=v_{j+1}v_{j+2}$ become red. Moreover, as $v_{j+2}$ is a white vertex in $G^D$, the edge $v_{j+2}v_{j+3}$ is either white or blue. In both cases, its weight decreases by at least $1$ when $e$ is played. Remark that $v_{j-1}v_{j}$ and $v_{j+2}v_{j+3}$  are different edges as the vertex $v_{j-1}$ is red while $v_{j+2}$ is white in $G^D$.
Then, the decrese in the weight can be estimated as $\w(e) \ge 1+ 2\cdot 2 +1=6$.  \smallqed
\end{itemize}

\begin{claim} \label{claim:4}
If $W\neq \emptyset$ and there are no adjacent white edges in $G^D$, then $\w(e)\ge 4$ for every playable edge $e$.  
\end{claim}
\proof Let $e \in W \cup B$. If $e$ is a white edge then, under the given conditions, every adjacent edge is from $B$. By Lemma~\ref{lem:2+}, there are two blue $C$-edges adjacent to $e$. Denote them by $e'$ and $e''$. By Claim~\ref{claim:1}, the edges $e$, $e'$, and $e''$ will be red in $G^{D\cup \{e\}}$. This shows $\w(e) \ge 2+2 \cdot 1=4$. Now suppose that $e \in B$  and it is adjacent to the white edge $f$ in $G^D$.  By Lemma~\ref{lem:2+}, there are two blue $C$-edges, say $f'$ and $f''$, which are adjacent to $f$. (Note that $e$ might be identical with $f'$ or $f''$.)  Again by Claim~\ref{claim:1}, the edges $f$, $f'$, and $f''$ will be red in $G^{D\cup \{e\}}$. This yields $\w(e) \ge 2+2 \cdot 1=4$. \smallqed
\medskip

Now, we are ready to finish the proof of Theorem~\ref{thm:edge-dom-game}. If the residual graph $G^D$ contains two adjacent white edges before Dominator's move $d_i$ then, by Claim~\ref{claim:3}, this move decreases the weight by at least $6$. If the game is not finished with $d_i$, according to Claim~\ref{claim:2}, Staller's move $s_i$ results in an additional decrease of at least $2$. If $G^D$ does not contain two adjacent white vertices then, by Claim~\ref{claim:4}, $\w(d_i)\ge 4$ and Staller's move $s_i$ further decreases the weight by at least $4$. Let $k$ denote the number of turns in the game. No matter whether Dominator or Staller finishes the game, the total decrease in the weight of $G$ is at least $4k$. On the other hand, the total decrease is $f(G^\emptyset)-0=2m(G)$. Thus, Dominator's greedy strategy ensures $4k \le 2m(G)$. As $\gamma_{e,g}(G) \le k$, we conclude $\gamma_{e,g}(G) \le m(G)/2$.  \qed

\subsection{Graphs with edge dominating trails}

In this subsection, we prove a theorem that is analogous to Theorem~\ref{thm:edge-dom-game} but concerns
 graphs that have an edge dominating open trail. In the proof, we refer to the edge domination game on a residual graph $F^D$ that is equivalently the game on the predominated graph $F|_{N[D]}$. The invariant $\gamma_{e,g}(F^D)$ denotes the number of edges played in an edge domination game starting with $F^D$, assuming that both players apply an optimal strategy and Dominator plays first.

\begin{theorem} \label{thm:2.9}
 If $G$ is a graph with an edge dominating trail, then $$\gamma_{e,g}(G)\le \left\lceil \frac{m(G)}{2}\right\rceil.$$
\end{theorem}	 
\proof If $G$ is a path on $n$ vertices, then $L(G)\cong P_{n-1}$ and, by Observation~\ref{obs:1}, we have $\gamma_{e,g}(G)=\gamma_{g}(P_{n-1}) \le \lceil(n-1)/2\rceil= \lceil m(G)/2 \rceil$. ((We have also used here the known formula for $\gamma_g(P_n)$, see~\cite{kosmrlj-2017}.) Note that if $G$ admits an edge dominating circuit, then Theorem~\ref{thm:2.9} is a direct consequence of Theorem~\ref{thm:edge-dom-game}.

From now on, assume that $G$ is not a path but has an edge dominating open trail $v_1 \dots v_\ell$. We first construct a graph $F$ from $G$ by adding two new vertices $v_1'$ and $v_\ell'$ and three edges $v_1v_1'$, $v_1' v_\ell'$,  $v_\ell v_\ell'$ to it. Note that $v_1 \dots v_\ell v_\ell' v_1'v_1$ is an edge dominating circuit in $F$. Furthermore, since $G$ is not a path, $F$ is not a cycle. By setting $D=\{v_1'v_\ell'\}$ we obtain the residual graph $F^D$. We first consider the edge domination game that starts on $F^D$ with a move of Dominator and prove the following estimation.
\begin{claim} \label{cl:2.10}
 $\gamma_{e,g}(F^D) \le \left\lceil \frac{m(G)}{2} \right\rceil$.
\end{claim}

\proof As $F^D$ is a residual graph obtained from $F$ which is not a cycle and has an edge dominating circuit, we may use the weight function $f$ from the proof of Theorem~\ref{thm:edge-dom-game} and apply Claims~\ref{claim:2}, \ref{claim:3} and \ref{claim:4}. Observe that  $f(F^D)=2m(G)+2$ and the weight of the residual graph equals $0$ at the end of the game. Let $k$ denote the number of played edges in the edge domination game on $F^D$. By Claims~\ref{claim:2}, \ref{claim:3} and \ref{claim:4},  Dominator has a startegy which ensures that the average decrease is at least $4$ in a move and therefore $f(F^D)-0 \ge 4k$. We conclude 
$$\gamma_{e,g}(F^D) \le k \le \left\lfloor \frac{2m(G)+2}{4} \right\rfloor =\left\lceil \frac{m(G)}{2} \right\rceil$$
that establishes the statement. \smallqed
\medskip

To finish the proof of the theorem, it sufficies to verify the following inequality.
\begin{claim} \label{cl:2.11}
	$\gamma_{e,g}(G) \le \gamma_{e,g}(F^D)$.
\end{claim}

\proof We apply the proof method named ``Imagination Strategy''. Let Game~1 be an edge domination game played on $G$ (i.e., on $G^\emptyset$) and Game~2  the imagined edge domination game started on $F^D$. Staller plays optimally in Game~1. After her each move in Game~1,  Dominator ``interprets'' it in Game~2. This means that if the same edge can be legally played in Game~2, then Dominator copies it there; otherwise Staller's move can be replaced by any legal move in Game~2. Then, Dominator replies optimally in Game~2 and also interprets it in Game~1. Again, this means that he copies the move into Game~1, if possible; otherwise Dominator plays an arbitrary legal move in Game~1. 

We may assume that Dominator never plays $v_1v_1'$ or $v_\ell v_\ell'$ in Game~2. Indeed, as $N_F[v_1v_1'] \cap W \subseteq N_F[v_1v_2]$, the move $v_1v_1'$ can be replaced by $v_1v_2$ in each optimal strategy of Dominator. The same is true for $v_\ell v_\ell'$ and $v_\ell v_{\ell-1}$. Under this assumption, we prove the following.
\begin{itemize}
	\item[($\star$)] After each move and its interpretation in the other game, the set $U_1$ of undominated edges in Game~1 is the same as the set $U_2$ of undominated edges in Game~2. 
\end{itemize}
Before the first move, we have $U_1=U_2= E(G)$. If $U_1 = U_2$ is true before a move of Dominator, then his optimal move (that is not $v_1v_1'$ or $v_\ell v_\ell'$) in Game~2 dominates an edge $e \in U_2$. Then, by our assumption,  $e \in U_1$ also holds and, consequently, the same move is legal in Game~1. Therefore, $U_1=U_2$ remains valid for the updated sets. Similarly, if $U_1=U_2$ holds before Staller's (optimal) move in Game~1, then the same edge is playable in Game~2 and after Dominator copies the move to Game~2, the two sets of undominated edges remain equal. This implies that Game~1 and Game~2 finishes at the same time, say after the $k^{\mbox{th}}$ move.

Since Staller plays optimally in Game~1 but Dominator might not, we have $\gamma_{e,g}(G) \le k$. Similarly, as Dominator plays optimally in Game~2 but Staller might not, we have $\gamma_{e,g}(F^D) \ge k$. The two inequalities together prove $\gamma_{e,g}(G) \le \gamma_{e,g}(F^D)$. \smallqed
\medskip

Our remarks on the case of $G \cong P_n$ together with Claims~\ref{cl:2.10} and \ref{cl:2.11} establish the statement  $\gamma_{e,g}(G)\le \left\lceil \frac{m(G)}{2}\right\rceil.$ \qed

\subsection{Traceable line graphs}

Based on the results in the previous subsections, we prove here that Rall's Conjecture is valid on the class of line graphs.

\begin{theorem} \label{thm:dom-game}
	If $G$ is a traceable line graph, then $$\gamma_{g}(G)\le \left\lceil \frac{n(G)}{2}\right\rceil.$$
\end{theorem}
\proof Since $G$ is a line graph, there exists an isolate-free graph $F$ such that $G=L(F)$. Clearly, we have $n(G)=m(F)$ and, by Observation~\ref{obs:1}, $\gamma_g(G)=\gamma_{e,g}(F)$. 
It is well-known (see e.g.\ {\rm \cite[Proposition 8]{HaNash}}) that if the line graph $G=L(F)$ is hamiltonian, then $F$ admits an edge dominating circuit. Similarly, if $G$ is traceable, then $F$ has an edge dominating trail. By Theorem~\ref{thm:2.9}, $\gamma_{e,g}(F) \le \lceil m(F)/2 \rceil$. Together with the previous observations it directly implies $\gamma_{g}(G) \le \lceil n(G)/2 \rceil$. \qed

\section{Claw-free graphs of minimum degree $2$} \label{sec:claw-free-mindeg-2}

In the previous sections, we proved that Conjecture~\ref{conj:1/2-for-all} is true for claw-free cubic graphs and also for Hamiltonian line graphs. Now, we consider a wider graph class, namely the class of claw-free graphs, and show that a weakened version of the conjecture remains true.

\begin{theorem} \label{thm:claw-free-mindeg-2}
	If $G$ is a claw-free graph and $\delta(G) \ge 2$, then $$\ggd(G) \le \left\lceil \frac{11}{20} \, n(G) \right\rceil.$$
 Moreover, if $G$ is also connected and it is neither a $5$-cycle nor a $9$-cycle, then 
	$$\ggd(G) \le \left\lfloor \frac{11}{20} \, n(G) \right\rfloor.$$
\end{theorem}

\proof 
We first recall that $\gamma_{g}(C_n)= \lceil\frac{n}{2} \rceil =\frac{n+1}{2}$ holds for a cycle $C_n$ if $n \equiv 1 \mod 4$, while $\gamma_{g}(C_n)= \lfloor\frac{n}{2} \rfloor$ is valid for the remaining cases.
Consequently, we have $\gamma_{g}(C_n) \le \frac{11}{20}n$ if $n \notin \{5,9\}$.
It is straightforward to check that $\gamma_{g}(C_n) \le \lceil \frac{11}{20}n \rceil$ holds for $n=5$ and $9$.  

For the rest of the proof set $n = n(G)$. If $G$ is a disjoint union of cycles, then applying~\cite[Corollary~18]{ruk-2019} we deduce that $\gamma_{g}(G) \le \lceil \frac{n}{2} \rceil \le \lceil \frac{11}{20}n \rceil $. From now on, we assume that $G$ satisfies the conditions of the theorem and it is not a disjoint union of cycles. We therefore have $\Delta(G) \ge 3$. Consider a D-game on $G$ and split it into two parts. Phase~$1$ starts with the first move $d_1$ of the game; Phase~$2$ begins with the first move of Dominator when the (vertex) residual graph before that move satisfies $\Delta_W(W)\le 2$.   Phase~$2$ ends when the game is finished. It might happen that the game ends in Phase~$1$ and no move belongs to Phase~$2$. In every residual graph $G^D$, no matter whether $G^D$ was obtained in Phase~$1$ or Phase~$2$, we call a blue vertex a B$^+$-vertex if it became blue in Phase~$1$ of the game, otherwise it is a B$^-$-vertex. $B^+$ and $B^-$ denote the corresponding sets of blue vertices in $G^D$. Clearly, $B=B^+$ if the residual graph belongs to Phase~$1$ of the game. Now, we define the potential function 
$$f(G^D)= 22|W|+10|B^+|+9|B^-|,$$ 
or equivalently, we assign a weight of $22$, $10$, $9$ and $0$, respectively, to each white, B$^+$-, B$^-$-, and red vertex of $G^D$. As it is defined in Section~\ref{sec:proof-meth},  $\sa(v)$ stands for the difference $f(G^D)-f(G^{D\cup\{v\}})$. We prove three claims, one for general consequences of claw-freeness, one for Staller's moves, and one for a possible strategy of Dominator.

\begin{claim} \label{claim:gen-5}
	Let $G^D$ be a residual graph that was obtained from a claw-free graph $G$ with $\delta(G) \ge 2$. 
	\begin{itemize}
	\item[$(a)$] If $v \in B$, then $N(v) \cap W$ induces a complete graph.
		\item[$(b)$] If $u$ and $v$ are adjacent vertices such that  $u \in W$, $v\in B$ and  $\deg_W(u) \ge 2$, then $v$ is a B$^+$-vertex.
		\item[$(c)$] If $u \in W$ and its only white neighbor is $u'$, then $u$ has a blue neighbor $v$ such that whenever $u$ becomes dominated in a later residual graph, $v$ becomes red at the same time.
	\end{itemize}	
\end{claim}
\proof $(a)$ Suppose to the contrary that $N(v) \cap W$ contains two independent vertices $u_1$ and $u_2$.  Since $v$ is blue, it is dominated by a vertex $x \in D$. Then, $x$ is red and $N(x)$ does not contain white vertices. Therefore, $u_1 \notin N(x)$ and $u_2 \notin N(x)$. We may conclude that the set $\{u_1, u_2, x\}$ is independent and hence that $v, u_1, u_2, x$ induce a claw in $G$. This contradiction proves $(a)$.

$(b)$ Consider the turn in the game when $v$ became blue. Then, $v$ was white in the residual graph $G^D$ just before this move. We infer that before the move $\deg_W(u)\ge 3$ and consequently, $\Delta_W(W) \ge 3$ was true in $G^D$. This implies that the move played in $G^D$ belongs to Phase~$1$ of the game and that $B^-=\emptyset$ remains valid when $v$ becomes blue. 

$(c)$ By definition, every neighbor of a white vertex belongs to $W \cup B$. Since $\deg_W(u) =1$ and $\deg(u) \ge 2$, there is at least one blue vertex $v$ in $N(u)$. By $(a)$, we know that $N(v) \cap W$ induces a complete graph. Since it also includes vertex $u$, there are only two possibilities, namely either $N(v) \cap W= \{u\} $ or $N(v) \cap W= \{u, u'\} $. In the former case, $u$ is the only white neighbor of $v$ and therefore, the statement of $(c)$ is valid for $v$. In the second case, $v$ does not satisfy the statement only if there exists a blue vertex $y$ that is adjacent with $u$ but not with $u'$. But in this case we have the setup $y \in B$, $N(y) \cap W= \{u\} $ for $y$ itself and so, $(c)$ is satisfied with $y$.
\smallqed

\begin{claim} \label{claim:st-5}
If Staller plays a vertex $v$ in the residual graph $G^D$, then $\sa(v) \ge 22$.
\end{claim}
\proof As $v$ is a legal move in $G^D$, it is either white or blue. If $v \in W$, then, as $v$ becomes red by this move, $f^D(v)$ decreases by $22$. This immediately implies $\sa(v) \ge 22$. If $v \in B^+$, then it dominates at least one white vertex, say $u$. Since $v$ becomes red and $u$ becomes blue or red,  $\sa(v) \ge (10-0)+(22-10)=22$. If $v \in B^-$, then Phase~$2$ has already begun and therefore, the white neighbor $u$ of $v$ cannot become a B$^+$-vertex. This gives $\sa(v) \ge (9-0)+(22-9)=22$ again. 
\smallqed

\begin{claim} \label{claim:dom-5}
	If $G^D$ is a residual graph in a D-game and Dominator is the next player with his i$^{\rm th}$ move, then the following statements are true.
 \begin{itemize} 
 \item[$(1)$] If $G^D[W]$ is not a disjoint union of cycles or if $B\neq \emptyset$, then
  \begin{itemize}
  \item[$(i)$] Dominator either can ensure that the potential function decreases by at least $80$ during the next two moves $d_i$ and $s_i$; or
  \item[$(ii)$]  Dominator finishes the game with the next move $d_i$ that decreases $f(G^D)$ by at least $40$.
  \end{itemize}
  \item[$(2)$] If $G^D[W]$ is a disjoint union of cycles and $B= \emptyset$, then Dominator can ensure that the game finishes in at most $\frac{|W|+1}{2}$ moves.
 \end{itemize}
\end{claim} 
\proof
$(2)$ The condition $B=\emptyset$ in statement $(2)$ implies that a D-game on $G^D$ is the same as a D-game on $G^D[W]$. Since the latter is a disjoint union of cycles on $|W|$ vertices, Dominator can ensure (applying \cite[Corollary~18]{ruk-2019} again) that the game ends in at most $\left\lceil \frac{|W|}{2}\right\rceil$ moves. This verifies the upper bound $\frac{|W|+1}{2}$ on the number of remaining moves.
\medskip

$(1)$ Let the residual graph $G^D$ satisfy the conditions given in $(1)$. We consider several cases in the proof.
\paragraph{(D1)}  $\Delta_W(W) \ge 3$. \\
Let $u$ be a white vertex that has at least three white neighbors. If Dominator plays $u$, then $u$ becomes red and the three neighbors become blue or red. In this way, Dominator can ensure that the decrease in the potential function is at least $\sa(u)\ge (22-0)+3(22-10)=58$. If this is the last move in the game, then $(ii)$ is satisfied. Otherwise, by Claim~\ref{claim:st-5}, Staller's next move decreases the potential function by at least $22$ and the two moves together cause a decrease of at least $80$ as stated.

\paragraph{(D2)} $\Delta_W(W) =2$ and $G^D[W]$ contains a component $P_\ell$ that is a path on at least $3$ vertices.\\
Suppose that $\{u_1,u_2,u_3\} \subseteq W$, $\deg_W(u_1)=1$ and that these three vertices induce a path. If Dominator plays $u_2$, then both $u_1$ and $u_2$ become red and $u_3$ becomes blue or red. These changes decrease the sum of the weights by at least $2\cdot 22+(22-10)=56$. Further, by Claim~\ref{claim:gen-5} $(c)$, there exists a blue neighbor $v$ of $u_1$ which turns red when $u_1$ is dominated. This decreases the weight of $v$ by at least $9$. We infer that $\sa(u_2) \ge 56+9=65$ and either $(i)$ or $(ii)$ follows as in the previous case.

\paragraph{(D3)} $\Delta_W(W) \le 2$ and there exist two vertices $u_1$ and $u_2$ which induce a $P_2$-component in $G^D[W]$ and $|(N(u_1) \cup N(u_2))\cap B| \ge 2$.\\
If Dominator plays $u_1$, then both $u_1$ and $u_2$ become red. Moreover, by Claim~\ref{claim:gen-5} $(a)$, for each blue neighbor $v$ of $u_i$, $i \in [2]$, the vertices in $N(v)\cap W$ form a complete graph. As follows, $N(v) \subseteq \{u_1, u_2\}$ and after playing $u_1$, the blue vertex $v$ becomes red. As it is also true for the other blue neighbor, we infer $\sa(u_1) \ge 2 \cdot 22+2 \cdot 9=62$. Together with Claim~\ref{claim:st-5} this fact implies the statement. 
 
\paragraph{(D4)} $\Delta_W(W) =2$ and there is a cycle component $C_\ell$ in $G^D[W]$ such that at least one vertex from $C_\ell$ is adjacent to a blue vertex.\\
Without loss of generality, we assume that $u_1$, $u_2$ and $u_3$ are three consecutive vertices from the cycle and $v$ is a blue neighbor of $u_2$. By Claim~\ref{claim:gen-5} $(a)$, all white neighbors of $v$ are from $\{u_1,u_2,u_3\}$. Thus, after playing $u_2$, the vertices $u_2$ and $v$ become red, while $u_1$ and $u_3$ become blue or red. Note that, since $\Delta_W(W) =2$, this move of Dominator belongs to Phase~$2$ of the game and no new B$^+$-vertices may arise. Observe further that, by Claim~\ref{claim:gen-5} $(b)$, $v$ is a B$^+$-vertex. We conclude $\sa(u_2) \ge 22+10+2\cdot(22-9)=58$ that verifies the statement again.

\paragraph{(D5)} If none of the conditions (D1)--(D4) holds, but it is still true that $G^D[W]$ is not a disjoint union of cycles or  $B\neq \emptyset$, we get that at least one component of $G^D[W]$ is not a cycle. As $\Delta_W(W) \le 2$, this component is a path. Excluding the cases (D2) and (D3), this component is either a single vertex $x$ forming a $P_1$ component or two adjacent vertices $u_1u_2$ that have exactly one (common) blue neighbor $v$. In the latter case, each (present or later) move that dominates at least one of $u_1$ and $u_2$ makes all the three vertices red, and decreases the weight of the residual graph by at least $2 \cdot 22+ 9=53$. In the former case, as $\deg_G(x)\ge 2$ and $\deg_W(x)=0$, the white vertex $x$ has at least two blue neighbors, say $v_1$ and $v_2$. Referring again to Claim~\ref{claim:gen-5} $(a)$,  we get that $N(v_i)\cap W=\{x\}$. Thus, each (present or later) move that dominates $x$, turns $x$, $v_1$ and $v_2$ red and decreases the weight of the residual graph by at least $22+2 \cdot 9=40$.

Under the conditions of (D5), Dominator plays a vertex from a non-cycle component of $G^D[W]$ and decreases $f(G^D)$ by at least $40$. If the game is over, then $(ii)$ holds. Otherwise, in the next turn, Staller either plays a vertex that dominates a $P_1$- or a $P_2$-component from $G^D[W]$ and decreases $f(G^D)$ by at least $40$, or she plays a vertex from a cycle component of $G^D[W]$. Recall that, by the exclusion of the case (D4), she does not have a possibility for playing an outer blue vertex and dominating only one or two vertices from the cycle. Therefore, in the latter case, Staller dominates three white vertices such that one of them becomes red. This move decreases $f(G^D)$ by at least $22+2\cdot (22-9)=48$. Under any legal choice of Staller, the two consecutive moves of the game decrease the potential function by at least $80$. This finishes the proof of $(1)$.
\smallqed
\medskip

To complete the proof of Theorem~\ref{thm:claw-free-mindeg-2}, we first assume that all moves of Dominator can be done under the condition of $(1)$ in Claim~\ref{claim:dom-5}. Observe that this assumption covers all cases when $G$ is connected and not a cycle. According to Claim~\ref{claim:dom-5} $(1)$, Dominator can ensure that if $k$ vertices are played in the game, then the original weight $f(G^\emptyset)=22n$ decreases by at least $40k$ until the end of the game. Since the game is over when $f(G^D)=0$, we conclude that $22n-0 \ge 40k$ and hence,
$$\gamma_{g}(G) \le k \le \left\lfloor \frac{22n}{40}
\right\rfloor = \left\lfloor \frac{11}{20}n \right\rfloor.$$

In the other case, there is a smallest index $(i+1)$ such that  Dominator plays his $(i+1)^{\rm st}$ move in a residual graph $G^D$ which satisfies $B=\emptyset$ and $G^D[W]$ consists of cycle components. Let us suppose that $k$ is the total number of vertices played in the game. Claim~\ref{claim:dom-5} $(1)$ says that the potential function decreased by at least $2i \cdot 40=80i$ during the first $2i$ turns. By Claim~\ref{claim:dom-5} $(2)$, all  $|W|$ vertices of the cycle components become red in at most $\frac{|W|+1}{2}$ moves. By introducing the notation $j=k-2i$, we get $j \le \frac{|W|+1}{2}$. This implies $2j-1 \le |W| $. Thus, during the last $j$ turns, $f(G^D)$ decreases by at least $22|W|\ge 44j-22$. 
Comparing the number of played vertices, that is $k =2i+j$, and the total decrease in the weights, we obtain
$$22n \ge 80i+44j-22=40k+4j-22 \ge 40k-18,$$
as $j \ge 1$. From this inequality we obtain the desired estimation for the game domination number of $G$ that finishes the proof of the theorem:
$$ \gamma_{g}(G) \le k \le \left\lfloor \frac{11n+9}{20} \right\rfloor \le \left\lceil \frac{11}{20}n \right\rceil.$$
\qed

\section{Concluding remarks}

Using a computer search we have checked that Conjecture~\ref{conj:1/2-for-all} is supported by all graphs $G$ on at most $10$ vertices. The computer search becomes slow for larger graphs, even when  restricted to connected graphs $G$ with diameter at least $3$, minimum degree smaller than $5$, and domination number such that $\lceil n(G)/2 \rceil < 2 \gamma(G) - 1$. However, we did try some random approaches and they have all supported the conjecture. 

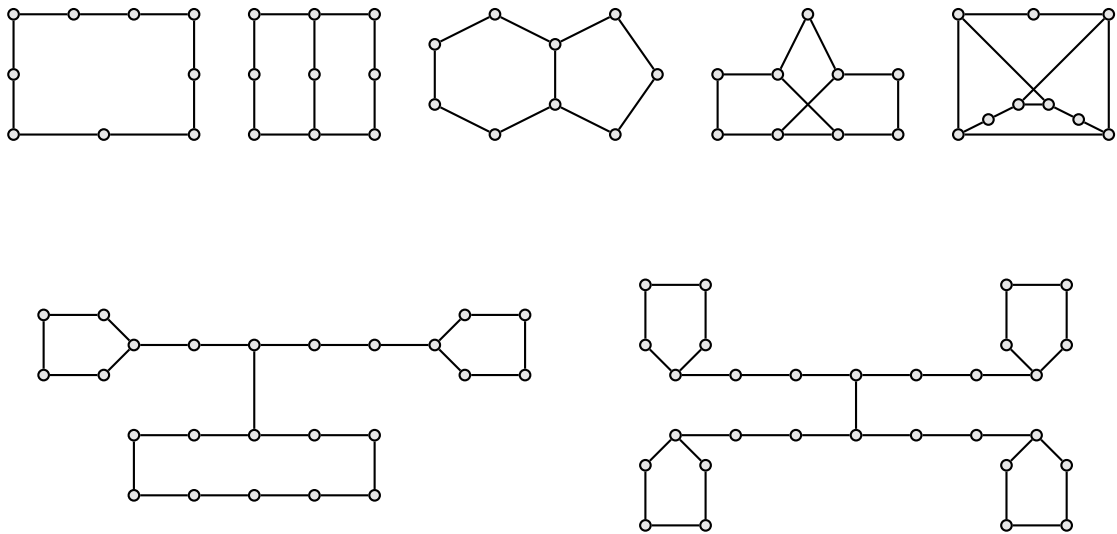
\begin{figure}[ht!]
	\begin{center}
		\begin{tikzpicture}[thick,scale=0.8]
		
		\tikzstyle{every node}=[circle, draw, fill=black!10,
		inner sep=0pt, minimum width=4pt]
		
		\begin{scope}
		\node (0) at (0,0) {};
		\node (1) at (1.5,0) {};
		\node (2) at (3,0) {};
		\node (3) at (3,1) {};
		\node (4) at (3,2) {};
		\node (5) at (2,2) {};
		\node (6) at (1,2) {};
		\node (7) at (0,2) {};
		\node (8) at (0,1) {};
		
		\draw (0) -- (1) -- (2) -- (3) -- (4) -- (5) -- (6) -- (7) -- (8) -- (0);
		\end{scope}
		
		\begin{scope}[xshift=4cm]
		\node (0) at (0,0) {};
		\node (1) at (1,0) {};
		\node (2) at (2,0) {};
		\node (3) at (2,1) {};
		\node (4) at (2,2) {};
		\node (5) at (1,2) {};
		\node (6) at (0,2) {};
		\node (7) at (0,1) {};
		\node (8) at (1,1) {};
		
		\draw (0) -- (1) -- (2) -- (3) -- (4) -- (5) -- (6) -- (7) -- (0);
		\draw (1) -- (8) -- (5);
		\end{scope}
		
		\begin{scope}[xshift=7cm]
		\node (0) at (0,0.5) {};
		\node (1) at (1,0) {};
		\node (2) at (2,0.5) {};
		\node (3) at (3,0) {};
		\node (4) at (3.7,1) {};
		\node (5) at (3,2) {};
		\node (6) at (2,1.5) {};
		\node (7) at (1,2) {};
		\node (8) at (0,1.5) {};
		
		\draw (0) -- (1) -- (2) -- (3) -- (4) -- (5) -- (6) -- (7) -- (8) -- (0);
		\draw (2) -- (6);
		\end{scope}
		
		\begin{scope}[xshift=11.7cm]
		\node (0) at (0,0) {};
		\node (1) at (1,0) {};
		\node (2) at (2,0) {};
		\node (3) at (3,0) {};
		\node (4) at (3,1) {};
		\node (5) at (2,1) {};
		\node (6) at (1.5,2) {};
		\node (7) at (1,1) {};
		\node (8) at (0,1) {};
		
		\draw (0) -- (1) -- (2) -- (3) -- (4) -- (5) -- (6) -- (7) -- (8) -- (0);
		\draw (1) -- (5);
		\draw (2) -- (7);
		\end{scope}
		
		\begin{scope}[xshift=15.7cm]
		\node (0) at (0,0) {};
		\node (1) at (2.5,0) {};
		\node (2) at (2.5,2) {};
		\node (3) at (1.25,2) {};
		\node (4) at (0,2) {};
		\node (5) at (0.5,0.25) {};
		\node (6) at (1,0.5) {};
		\node (7) at (1.5,0.5) {};
		\node (8) at (2,0.25) {};
		
		\draw (0) -- (1) -- (2) -- (3) -- (4) -- (0);
		\draw (0) -- (5) -- (6) -- (2);
		\draw (1) -- (8) -- (7) -- (4);
		\draw (6) -- (7);
		\end{scope}
		
		\begin{scope}[xshift=2cm, yshift=-6cm]
		\node (0) at (0,0) {};
		\node (1) at (1,0) {};
		\node (2) at (2,0) {};
		\node (3) at (3,0) {};
		\node (4) at (4,0) {};
		\node (9) at (0,1) {};
		\node (8) at (1,1) {};
		\node (7) at (2,1) {};
		\node (6) at (3,1) {};
		\node (5) at (4,1) {};
		\node (10) at (2,2.5) {};
		\node (11) at (1,2.5) {};
		\node (12) at (0,2.5) {};
		\node (13) at (-0.5,2) {};
		\node (14) at (-1.5,2) {};
		\node (15) at (-1.5,3) {};
		\node (16) at (-0.5,3) {};
		\node (17) at (3,2.5) {};
		\node (18) at (4,2.5) {};
		\node (19) at (5,2.5) {};
		\node (20) at (5.5,2) {};
		\node (21) at (6.5,2) {};
		\node (22) at (6.5,3) {};
		\node (23) at (5.5,3) {};
		
		\draw (0) -- (1) -- (2) -- (3) -- (4) -- (5) -- (6) -- (7) -- (8) -- (9) -- (0);
		\draw (7) -- (10) -- (11) -- (12) -- (13) -- (14) -- (15) -- (16) -- (12);
		\draw (10) -- (17) -- (18) -- (19) -- (20) -- (21) -- (22) -- (23) -- (19);
		\end{scope}
		
		\begin{scope}[xshift=11cm, yshift=-5cm]
		\node (0) at (0,0) {};
		\node (1) at (1,0) {};
		\node (2) at (2,0) {};
		\node (3) at (3,0) {};
		\node (4) at (4,0) {};
		\node (5) at (5,0) {};
		\node (6) at (6,0) {};
		\node (7) at (0,1) {};
		\node (8) at (1,1) {};
		\node (9) at (2,1) {};
		\node (10) at (3,1) {};
		\node (11) at (4,1) {};
		\node (12) at (5,1) {};
		\node (13) at (6,1) {};
		
		\node (14) at (-0.5,-0.5) {};
		\node (15) at (-0.5,-1.5) {};
		\node (16) at (0.5,-1.5) {};
		\node (17) at (0.5,-0.5) {};
		
		\node (18) at (5.5,-0.5) {};
		\node (19) at (5.5,-1.5) {};
		\node (20) at (6.5,-1.5) {};
		\node (21) at (6.5,-0.5) {};
		
		\node (22) at (-0.5,1.5) {};
		\node (23) at (-0.5,2.5) {};
		\node (24) at (0.5,2.5) {};
		\node (25) at (0.5,1.5) {};
		
		\node (26) at (5.5,1.5) {};
		\node (27) at (5.5,2.5) {};
		\node (28) at (6.5,2.5) {};
		\node (29) at (6.5,1.5) {};
		
		\draw (0) -- (1) -- (2) -- (3) -- (4) -- (5) -- (6);
		\draw (7) -- (8) -- (9) -- (10) -- (11) -- (12) -- (13);
		\draw (3) -- (10);
		
		\draw (0) -- (14) -- (15) -- (16) -- (17) -- (0);
		\draw (6) -- (18) -- (19) -- (20) -- (21) -- (6);
		\draw (7) -- (22) -- (23) -- (24) -- (25) -- (7);
		\draw (13) -- (26) -- (27) -- (28) -- (29) -- (13);
		
		\end{scope}
		
		\end{tikzpicture}
		\caption{The first row contains graphs $G$ with $n(G) = 9$ and $\ggd(G) = 5$. The second row contains a graph on $24$ vertices with $\ggd = 12$, and a graph on $30$ vertices with $\ggd = 15$.}
		\label{fig:eq}
	\end{center}
\end{figure}

If Conjecture~\ref{conj:1/2-for-all} is true, then it is best possible as demonstrated by cycles. Moreover, among diameter $2$ graphs, there are exactly eight equality cases, see~\cite[Fig.~1]{bujtas-2020c}. In addition, we were also able to find equality cases with diameter at least $3$. Using a computer search among graphs on at most $9$ vertices (and diameter at least $3$), we found exactly $5$ equality cases on $6$ vertices, $23$ equality cases on $8$ vertices, and $5$ equality cases on $9$ vertices. The latter are also the smallest (except for $C_5$) equality cases with odd number of vertices and are presented in the first row of Fig.~\ref{fig:eq}. Note that all these graphs are traceable. To add non-traceable equality cases, we present two sporadic examples, see the second row of Fig.~\ref{fig:eq}.

To conclude the paper we present two alternative versions of Conjecture~\ref{conj:1/2-for-all} that could become interesting provided a counterexample to Conjecture~\ref{conj:1/2-for-all} will be found.  

\begin{conjecture}
	If $\delta(G) \ge 2$, then $\gamma_g(G) \le  \frac{n(G)}{2} + C$, where $C$ is a universal constant. 
\end{conjecture}

From the examples above we know that $C \ge \frac{1}{2}$. A further weakening of the conjecture is the following. 

\begin{conjecture}
\label{con:constant-c}
	There exists a constant $c < \frac{3}{5} $ such that every graph $G$ with $\delta(G) \ge 2$ and $n(G) \ge 6$ satisfies $\gamma_g(G) \le c\cdot n(G)$.
\end{conjecture}

Recall that if $\delta(G) \ge 3$, then $\gamma_g(G) \le 0.5574\cdot n(G)$, see~\cite{bujtas-2015b}. Hence Conjecture~\ref{con:constant-c} holds for graphs of minimum degree at least $3$ (with the constant $c=0.5574$). 

\section*{Acknowledgements}

We are grateful to Ga\v{s}per Ko\v{s}mrlj for providing us with his software that computes game domination invariants. We acknowledge the financial support from the Slovenian Research Agency (research core funding No.\ P1-0297 and projects J1-9109, J1-1693, N1-0095, N1-0108).

\end{document}